\documentclass[12pt,a4paper]{amsart} 

\usepackage[T1]{fontenc} 
\usepackage[latin1]{inputenc} 
\usepackage{amsmath} 
\usepackage{amsthm} 
\usepackage{amsfonts} 
\usepackage{bbm} 
\usepackage{varioref} 
\usepackage{enumerate} 

\usepackage{ae}

\theoremstyle{definition}

\theoremstyle{remark}

\theoremstyle{plain}
\newtheorem{thm}{Theorem}
\newtheorem{lem}[thm]{Lemma}

\newcommand{\norm}[1]{\ensuremath{\left\Vert #1 \right\Vert}}
\newcommand{\abs}[1]{\ensuremath{\left\vert #1 \right\vert}}

\newcommand{\infabs}[1]{\ensuremath{\left\vert #1 \right\vert_\infty}}
\newcommand{\natnum}{\ensuremath{\mathbb{N}}}

\DeclareMathOperator{\dist}{dist}

\newcommand{\hdim}{\text{\textup{dim}}_{\text{\textup{H}}}}

\DeclareMathOperator{\Mat}{Mat}

\title[Restricted Diophantine approximation]{A metric theorem for
  restricted Diophantine approximation in positive characteristic}
\author[Simon Kristensen]{Simon Kristensen (Edinburgh)}
\date{\today}
\thanks{The author is a William Gordon Seggie Brown Fellow}

\address{Simon Kristensen, School of Mathematics, The University of
  Edinburgh, James Clerk Maxwell Building, Kings Buildings, Mayfield
  Road, Edinburgh, EH9 3JZ, Scotland}
\email{Simon.Kristensen@ed.ac.uk}

\dedicatory{To Maurice Dodson on his retirement}

\begin{document}

\maketitle

\section{Introduction}
\label{sec:introduction}

Let $\mathbb{F}$ be the finite field with $k = p^l$ elements, let
$\mathbb{F}[X]$ denote the ring of polynomials with coefficients from
$\mathbb{F}$, let $\mathbb{F}(X)$ denote the field of fractions over
this ring and let $\mathbb{F}((X^{-1}))$ denote the field of formal
Laurent series with coefficients from $\mathbb{F}$, \emph{i.e.},
\begin{equation*}
  \mathbb{F}((X^{-1})) = \left\{\sum_{i=n}^\infty a_i X^{-i} : n \in
  \mathbb{Z}, a_i \in \mathbb{F}\right\}.
\end{equation*}
We may define a non-Archimedean absolute value on $\mathbb{F}((X^{-1}))$
by 
\begin{equation*}
  \abs{\sum_{i=-n}^\infty a_i X^{-i}} = 
  \begin{cases}
    0 & \text{whenever } a_i = 0 \text{ for all } i \in \mathbb{Z}, \\
    k^{n} & \text{whenever } a_n \neq 0 \text{ and } a_i = 0 \text{
      for } i < n.
  \end{cases}
\end{equation*}
We can interpret $\mathbb{F}((X^{-1}))$ as the completion of
$\mathbb{F}(X)$ with respect to this absolute value. Note that in
addition to the usual non-Archimedean property of the absolute value,
$\abs{f+g} \leq \max\{\abs{f}, \abs{g}\}$, we also have
\begin{equation}
  \label{eq:27}
  \abs{f} \neq \abs{g} \Rightarrow \abs{f+g} = \max\{\abs{f},
  \abs{g}\}. 
\end{equation}

Diophantine approximation in $\mathbb{F}((X^{-1}))$, where a generic
element is approximated by elements from the field of fractions
$\mathbb{F}(X)$, has been studied by numerous authors (the survey
papers \cite{MR2001k:11135, MR2001j:11063} contain some of the known
results). Broadly speaking, the object of study has been variations
over inequalities of the form
\begin{equation*}
  \abs{f - \dfrac{P}{Q}} < \psi(\abs{Q}),
\end{equation*}
where $f \in \mathbb{F}((X^{-1}))$ and $P, Q \in \mathbb{F}[X]$ with
$Q \neq 0$. The study of the metric theory of Diophantine
approximation in this setting, in which the measure and Hausdorff
dimension of sets arising from such inequalities is studied, was begun
by de Mathan in \cite{MR43:161}, who proved an analogue of
Khintchine's theorem in Diophantine approximation. The author extended
this theorem to systems of linear forms \cite{kristensen03}.

Let $m,n \in \mathbb{N}$ and $\psi : \mathbb{F}[X]^m \rightarrow \{k^r
: r \in \mathbb{Z}\}$ be some function. Let $S \subseteq
\mathbb{F}[X]^m$.  For $\mathbf{v} \in \mathbb{F}((X^{-1}))^n$, let
$\norm{\mathbf{v}}$ denote the distance from $\mathbf{v}$ to the
nearest element in $\mathbb{F}[X]^n$ with respect to metric induced by
the norm $\infabs{\mathbf{x}} = \max\{\abs{x_1}, \dots, \abs{x_n}\}$,
where $\mathbf{x} = (x_1, \dots, x_n) \in \mathbb{F}((X^{-1}))^n$. In
this paper, we study the set
\begin{multline}
  \label{eq:1}
  \mathcal{W}_S(m,n;\psi) = \big\{A \in \Mat_{m \times
    n}\left(\mathbb{F}((X^{-1}))\right) : \norm{\mathbf{q} A}
  <  \psi(\mathbf{q}), \\
  \text{ for infinitely many } \mathbf{q} \in S \big\},
\end{multline}
where $\Mat_{m \times n}\left(\mathbb{F}((X^{-1}))\right)$ denotes the
space of $m$ by $n$ matrices $A$ with coefficients from
$\mathbb{F}((X^{-1}))$ and $\mathbf{q} A$ denotes the usual matrix
product.

Recently, Inoue and Nakada proved a Khintchine type theorem for the
special case $\mathcal{W}_Q(1,1;\psi)$ \cite[Theorem
1]{inoue03:_dioph}. Namely, they showed that the Haar measure of
$\mathcal{W}_Q(1,1;\psi)$ is null or full accordingly as the series
$\sum_{q \in S} \abs{q} \psi(q)$ converges or diverges, under the
additional assumptions that the approximating fractions $P/Q$ are on
lowest terms and that $\psi(q)$ depends only on $\abs{q}$.

In the real case, the Hausdorff dimension of the analogous sets for $m
= 1$ and arbitrary $n$ in the special case when $\psi = \abs{q}^{-v}$
was determined by Borosh and Fraenkel \cite{MR46:7175}. Various more
general cases were studied using the notion of ubiquitous systems by
Rynne \cite{MR93a:11066}. This was subsequently generalised to an even
more general form of approximation by Dickinson and Rynne
\cite{MR2001g:11129}.

We will consider the analogue of the case originally considered by
Rynne \cite{MR93a:11066}. It is the purpose of the present paper to
determine the Hausdorff dimension -- and in some cases the Haar
measure -- of the sets $\mathcal{W}_S(m,n;\psi)$, subject to very mild
constraints on the approximation function $\psi$. We will denote by
$\mu$ the Haar measure on $\mathbb{F}((X^{-1}))$, normalised so that
the unit ball,
\begin{equation*}
  U = \left\{f \in \mathbb{F}((X^{-1})) : \abs{f} < 1\right\},
  \subseteq \mathbb{F}((X^{-1})),
\end{equation*}
has measure $1$. By abuse of notation, the normalised Haar measure on
vector spaces $V$ over $\mathbb{F}((X^{-1}))$ will also be denoted by
$\mu$. Similarly, we will denote by $U$ the unit cube in $V$,
\emph{i.e.}, the $\dim(V)$-fold Cartesian product of $U$ with itself.
This should cause no confusion.

We need a few definitions. First, we need an appropriate notion of the
exponent of convergence for the sequence $S$ given by
\begin{equation}
  \label{eq:2}
  v(S) = \inf \left\{v \in \mathbb{R} : \sum_{\mathbf{q} \in S}
    \infabs{\mathbf{q}}^{-v} < \infty \right\}.
\end{equation}
We also need the appropriate notion of the order at infinity
$\lambda(\psi)$ of the error function $\psi$,
\begin{equation}
  \label{eq:3}
  \lambda(\psi) = \lim_{\substack{\infabs{\mathbf{q}} \rightarrow
      \infty \\ \mathbf{q} \in S}} \dfrac{-\log
    \psi(\mathbf{q})}{\log \infabs{\mathbf{q}}},
\end{equation}
defined whenever it exists. We can now state the first main theorem.
\begin{thm}
  \label{thm:Main}
  Let $\psi : \mathbb{F}[X]^m \rightarrow \{k^r : r \in \mathbb{Z}\}$.
  Suppose that $\psi(\mathbf{q})$ depends only on
  $\infabs{\mathbf{q}}$, is decreasing as a function of
  $\infabs{\mathbf{q}}$ and that the order at infinity of the
  function $\psi$ exists.
  \begin{enumerate}
  \item \label{item:1} If $n \lambda(\psi) < v(S)$, then
    $\mathcal{W}_S(m,n;\psi)$ is full with respect to the Haar measure
    on $\Mat_{m \times n}(\mathbb{F}((X^{-1})))$.
  \item \label{item:2}
    If $n \lambda(\psi) \geq v(S)$, then 
    \begin{equation*}
      \hdim(\mathcal{W}_S(m,n;\psi)) =
      n(m-1) + \dfrac{n + v(S)}{1 + \lambda(\psi)},
    \end{equation*}
    where $\hdim(E)$ denotes the Hausdorff dimension of the set $E$. 
  \end{enumerate}
\end{thm}
Note that while we can calculate the dimension for each value of
$v(S)$ and $\lambda(\psi)$, we are not able to show in general
whether the measure of $\mathcal{W}_S(m,n;\psi)$ is null or full for
the critical value $n \lambda(\psi)= v(S)$ with the methods of the
present paper.

In analogy with a result of Rynne \cite{MR1633797}, we may deduce from
Theorem \ref{thm:Main} the Hausdorff dimension of the set
$\mathcal{W}_S(m,n;\psi)$ with very mild conditions on $\psi$. Indeed,
it suffices to assume that $\psi$ is bounded, non-negative and that
$\psi(\mathbf{q}) > 0$ for infinitely many $\mathbf{q}$. Note that in
this case by defining
\begin{equation*}
  \hat{\psi}(\mathbf{q}) = 
  \begin{cases}
    \psi(\mathbf{q}) & \text{whenever } \mathbf{q} \in S \\
    0 & \text{otherwise,}
  \end{cases}
\end{equation*}
we would have $\mathcal{W}_S(m,n;\psi) =
\mathcal{W}_{\mathbb{F}[X]^m}(m,n; \hat{\psi})$. Consequently, for
such more general $\psi$, we omit the set $S$ from our notation and
talk about $\mathcal{W}(m,n;\psi)$. We define
\begin{equation}
  \label{eq:14}
  \eta(\psi) = \inf\left\{\eta \in \mathbb{R} : \sum_{\mathbf{q} \in
      \mathbb{F}[X]^m} \infabs{\mathbf{q}}^n \left(
      \dfrac{\psi(\mathbf{q})}{\infabs{\mathbf{q}}} \right)^\eta <
    \infty \right\}. 
\end{equation}

\begin{thm}
  \label{thm:Main2}
  Suppose that $\psi: \mathbb{F}[X]^m \rightarrow \{k^r : r \in
  \mathbb{Z}\} \cup \{0\}$ is bounded, non-negative and that
  $\psi(\mathbf{q}) > 0$ for infinitely many $\mathbf{q}$. Then,
  \begin{equation*}
    \hdim(\mathcal{W}(m,n;\psi)) = n(m-1) + \min\{\eta(\psi), n\}.
  \end{equation*}
\end{thm}

Theorem \ref{thm:Main2} generalises part (\ref{item:2}) of Theorem
\ref{thm:Main}.  As before, an optimal condition for when the measure
is full does not follow from our approach. This is however a very
difficult problem generalising the Duffin--Schaeffer conjecture (see
\emph{e.g.} \cite{MR548467}) to systems of linear forms over
$\mathbb{F}((X^{-1}))$. For simultaneous approximation, such a result
is known for reals \cite{MR1099767} and has been announced for formal
power series \cite{MR2019008}.

\section{Proof of Theorem \ref{thm:Main}}
\label{sec:proof-theorem}

We need to prove three things. First, we will show that the right hand
side in \eqref{item:2} is an upper bound on the Hausdorff dimension of
$\mathcal{W}_S(m,n;\psi)$. Note that this implies that the Haar
measure is zero when $n \lambda(\psi) > v(S)$. Subsequently, we need
to show that the measure is full in case \eqref{item:1} and that the
right hand side in \eqref{item:2} is also a lower bound on the
dimension.

In the following, matrices in $\Mat_{m \times
  n}(\mathbb{F}((X^{-1})))$ will be identified with vectors in
$\mathbb{F}((X^{-1}))^{mn}$. Given a vector $\mathbf{x}$ and a set
$V$, $\dist(\mathbf{x},V)$ will denote the minimal distance from
$\mathbf{x}$ to $V$ in the absolute value $\infabs{\cdot}$. Given two
real numbers $a,b$, we will use the Vinogradov notation and say that
$a \ll b$ if there is a constant $K > 0$ such that $a \leq K b$. If $a
\ll b$ and $b \ll a$, we will write $a \asymp b$.

\subsection{An upper bound}
\label{sec:an-upper-bound}

We note that the set $\mathcal{W}_S(m,n;\psi)$ is
invariant under translations by elements from $\Mat_{m \times
  n}\left(\mathbb{F}[X]\right)$. Hence, we restrict ourselves to
considering the intersection of $\mathcal{W}_S(m,n;\psi)$
and the unit cube $U$. We will prove that the upper bound is the right
one when $n = 1$. In this case, we are determining the dimension of
the set
\begin{equation*}
  \mathcal{W}^*_S(m,1;\psi) = \big\{A \in U :
  \norm{\mathbf{q} A} < \psi(\infabs{\mathbf{q}})
  \text{ for infinitely many } \mathbf{q} \in S \big\}.
\end{equation*}
We omit the details of the case $n > 1$ for ease of notation, but give
an outline of the differences from the one-dimensional case at the end
of this part of the proof.

Consider the $(m-1)$-dimensional affine subspaces for which the left
hand side is equal to zero, \emph{i.e.}, the affine subspaces
\begin{equation*}
  H(\mathbf{q},p) = \left\{A \in U : \mathbf{q} A = p \right\}.
\end{equation*}
Note that for this set to be non-empty, $\abs{p} \leq
\infabs{\mathbf{q}}$. Clearly, points satisfying the relevant
inequality for a fixed $p$ and $\mathbf{q}$ must lie within
$\psi(\infabs{\mathbf{q}}) \infabs{\mathbf{q}}^{-1}$ of
$H(\mathbf{q},p)$. 

Let $\epsilon > 0$ be arbitrary. By definition, for
$\infabs{\mathbf{q}}$ large enough, $\psi(\infabs{\mathbf{q}}) \leq
\infabs{\mathbf{q}}^{-\lambda(\psi)+\epsilon}$. Suppose that
$\infabs{\mathbf{q}}$ is large enough for this to hold.

We cover the neighbourhoods of $H(\mathbf{q}, p)$ by $\asymp
\infabs{\mathbf{q}}^{(1+\lambda(\psi)-\epsilon)(m-1)}$ balls of radius
$2 \infabs{\mathbf{q}}^{-\lambda(\psi)+\epsilon-1}$. Call this cover
$\mathcal{C}(\mathbf{q},p)$. For any $M > 0$, the sets
$\bigcup_{\infabs{\mathbf{q}} \geq M, \mathbf{q} \in S}
\bigcup_{\abs{p} \leq \infabs{\mathbf{q}}} \mathcal{C} (\mathbf{q},p)$
cover $\mathcal{W}^*_S(m,1;\psi)$. The $s$-length of this cover
$\mathcal{C}(M)$ is
\begin{equation}
  \label{eq:28}
  \begin{split}
    \ell^s\left(\mathcal{C} (M) \right) &\ll
    \sum_{\substack{\infabs{\mathbf{q}} \geq M \\ \mathbf{q} \in S}}
    \sum_{\abs{p} \leq \infabs{\mathbf{q}}}
    \infabs{\mathbf{q}}^{(1+\lambda(\psi)-\epsilon)(m-1)}
    \infabs{\mathbf{q}}^{-(\lambda(\psi)-\epsilon+1)s} \\
    &\ll \sum_{\substack{\infabs{\mathbf{q}} \geq M \\ \mathbf{q} \in
        S}} \infabs{\mathbf{q}}^{1 + (1+\lambda(\psi)-\epsilon)(m-1) 
      -(\lambda(\psi)-\epsilon+1)s}\\
    &\ll \sum_{\substack{\infabs{\mathbf{q}} \geq M \\ \mathbf{q} \in S}}
    \infabs{\mathbf{q}}^{-v(S) + \epsilon'},
  \end{split}
\end{equation}
whenever
\begin{equation}
  \label{eq:4}
  s > (m-1) + \dfrac{n + v(S)}{1 + \lambda(\psi)}, 
\end{equation}
for some $\epsilon' > 0$, which tends to zero as $\epsilon$ does.
Hence, we may choose $\epsilon > 0$ (which was arbitrary) such that
the last series of \eqref{eq:28} tends to zero as $M$ tends to
infinity. By the Hausdorff--Cantelli Lemma (see \emph{e.g.}
\cite{MR2001h:11091}), the Hausdorff dimension of
$\mathcal{W}^*_S(m,1;\psi)$ must then be less than or equal to the
right hand side of \eqref{eq:4}.

To prove the statement for $n > 1$, we note that the multidimensional
analogues of $H(\mathbf{q}, p)$ will be $n(m-1)$-dimensional
affine spaces. We may again cover the neighbourhoods by balls of
radius $\asymp \infabs{\mathbf{q}}^{-\lambda(\psi)+\epsilon-1}$.  Using
elementary upper bounds for the $s$-length of the resulting cover
implies the result.

\subsection{Reduction to simpler $\psi$}
\label{sec:reduct-simpl-boldsym}

We will now show that it suffices to consider the case when
$\psi(\infabs{\mathbf{q}})$ is of the form
$\infabs{\mathbf{q}}^{-\lambda(\psi)}$. For this, we use again the
existence of $\lambda(\psi)$. By \eqref{eq:3}, for any $\epsilon > 0$
there is an $r_0 \in S$ such that for any $\mathbf{q} \in S$ with
$\infabs{\mathbf{q}} \geq k^{r_0}$,
\begin{equation*}
  \infabs{\mathbf{q}}^{-\lambda(\psi)- \epsilon} \leq
  \psi(\mathbf{q}).
\end{equation*}
We for define $v>0$ the set
\begin{multline}
  \label{eq:6}
  \mathcal{W}_S(m,n;v) = \big\{A \in \Mat_{m \times
    n}\left(\mathbb{F}((X^{-1}))\right) : \norm{\mathbf{q} A}
  < \infabs{\mathbf{q}}^{-v}, \\
  \text{for infinitely many } \mathbf{q} \in S
  \big\},
\end{multline}
By the above, for any $\epsilon > 0$,
\begin{equation*}
  \mathcal{W}_S(m,n;\lambda(\psi_1) + \epsilon) \subseteq
  \mathcal{W}_S(m,n;\psi). 
\end{equation*}
As $\epsilon > 0$ is arbitrary, it therefore suffices to study the
sets $\mathcal{W}_S(m,n;v)$ and prove the corresponding full
measure result and lower bound on the Hausdorff dimension for this
set.

\subsection{Measures and counting lemmas}
\label{sec:meas-count-lemm}

We will use some probabilistic lemmas to prove the second and third
part of the theorem. The method which we will use is adapted from that
used by Dodson in \cite{MR95d:11092}. First, we need some estimates
for the measure and number of elements of various sets. We define for
$\mathbf{q} \in S$ and $\epsilon > 0$ the set
\begin{equation}
  \label{eq:9}
  B(\mathbf{q},\epsilon) = \left\{ A \in U : \norm{\mathbf{q} A}
    < \epsilon \right\}. 
\end{equation}
From \cite[equation (2.6)]{kristensen03} we extract the following
lemma: 
\begin{lem}
  \label{lem:independence}
  Let $\epsilon, \epsilon' > 0$ and let $m \geq 2$. Suppose that
  $\mathbf{q}, \mathbf{q}' \in \mathbf{F}[X]^m$ are linearly
  independent over $\mathbb{F}((X^{-1}))$. Then,
  \begin{equation*}
    \mu\left(B(\mathbf{q},\epsilon) \cap B(\mathbf{q}', \epsilon')
    \right) = \mu\left(B(\mathbf{q}, \epsilon) \right) \mu \left(
      B(\mathbf{q}', \epsilon')\right). 
  \end{equation*}
\end{lem}
Furthermore, from \cite[equation (2.3)]{kristensen03}, we get
\begin{lem}
  \label{lem:measure}
  For any $\mathbf{q} \in S$,
  \begin{equation*}
    \mu\left(B(\mathbf{q},\infabs{\mathbf{q}}^{-v})\right) \asymp
    \infabs{\mathbf{q}}^{-vn}.
  \end{equation*}
\end{lem}

In general, measuring the intersection of the sets defined in
\eqref{eq:9} is difficult. However, we may restrict ourselves to
certain subsets for which the task is easier. When $m = 1$, we define
\begin{equation*}
  B'(q,\epsilon) = \left\{ A \in U : \infabs{q A - \mathbf{p}}
    < \epsilon \quad \text{for $(q, p_i) = 1$ for all $1 \leq i \leq
      n$} \right\}.  
\end{equation*}
We will also need an analogue of the Euler totient function, defined
for any $q \in \mathbb{F}[X]$ to be
\begin{equation}
  \label{eq:17}
  \Phi(q) = \# \{q' \in \mathbb{F}[X] : \abs{q'} < \abs{q}, q'
  \text{monic and } (q,q')=1\}
\end{equation}
From \cite[equation (9)]{MR2019008}, we distill the following lemma:
\begin{lem}
  \label{lem:1D-measure}
  For $\epsilon > 0$ small enough, 
  \begin{equation*}
    \mu(B'(q,\epsilon)) = \epsilon^n \dfrac{\Phi(q)^n}{\abs{q}^n}.
  \end{equation*}
\end{lem}
Additionally, we distill from \cite[page 159]{MR2019008} that
\begin{lem}
  \label{lem:1D-intersection}
  For $q, q' \in \mathbb{F}[X]$ and $\epsilon, \epsilon' > 0$, 
  \begin{equation*}
    \mu\left(B'(q,\epsilon) \cap B'(q',\epsilon')\right) \ll
    \epsilon^n {\epsilon'}^n.
  \end{equation*}
\end{lem}

We will now estimate the growth of the function $\Phi(q)$.
\begin{lem}
  \label{lem:phi-function}
  \begin{equation*}
    \Phi(q) \asymp \abs{q}.
  \end{equation*}
\end{lem}

\begin{proof}
  The proof is simple. For the purposes of the present proof, $\mu$
  denotes the usual M\"obius function. We use the fact that $\sum_{d
    \vert r} \mu(d)$ is zero whenever $r \neq 1$ and is one when
  $r=1$.  Now, for $q \in \mathbb{F}[X]$ with $\deg(q) = n$,
  \begin{multline*}
    \Phi(q) = \sum_{\substack{\deg (p,q) = 0 \\ p \text{ monic} \\
        \abs{p} < \abs{q}}} 1 = \sum_{\substack{\deg (p,q) = r \\ p
        \text{ monic} \\ \abs{p} < \abs{q}}} \sum_{d \vert r+1}
    \mu(d)\\ 
    = \sum_{d=1}^{n} \mu(d) \sum_{\substack{\abs{p} = k^{n+1-d} \\ p 
        \text{ monic}}} 1 = \sum_{d=1}^{n} \mu(d) k^{n+1-d} \gg k^n
    = \abs{q}.
  \end{multline*} 
  The $\gg$ follows as the preceding expression is a polynomial with
  leading term $k^n$. The converse upper inequality is trivial.
\end{proof} 

We will need the higher dimensional analogue of the sets
$B'(\mathbf{q}, \epsilon)$ in the course of the proof.  We will impose
a further restriction on the resonant neighbourhoods in order to
enable us to treat this situation. In particular, we need to control
the measure of the overlaps when $\mathbf{q}$ and $\mathbf{q}'$ are
linearly dependent. For $\mathbf{q} \in \mathbb{F}[X]^m$ and $\epsilon
> 0$, we define
\begin{multline*}
  B''(\mathbf{q},\epsilon) = \big\{ A \in U : \infabs{\mathbf{q} A -
    \mathbf{p}} < \epsilon \\
  \text{for $(\gcd(q_1, \dots, q_m), p_i) = 1$ for all $1 \leq i
    \leq n$} \big\}.
\end{multline*}
Here, $\gcd(q_1, \dots, q_m)$ denotes the greatest common denominator
of the coordinates of $\mathbf{q}$ in $\mathbb{F}[X]$, which is unique
up to multiplication by an element of $\mathbb{F}$. For these sets,
we also estimate the relevant measures
\begin{lem}
  \label{lem:big_m_measure}
  For $\epsilon > 0$ small enough, 
  \begin{equation*}
    \mu(B''(\mathbf{q},\epsilon)) \asymp \epsilon^n.
  \end{equation*}
\end{lem}

\begin{proof}
  This is shown exactly as Lemma \ref{lem:1D-measure}, by measuring
  each component of the set and summing over them. The asymptotic
  estimate follows on using Lemma \ref{lem:phi-function}.
\end{proof}

\begin{lem}
  \label{lem:big_m_intersection}
  Let $\mathbf{q}, \mathbf{q}' \in \mathbb{F}[X]^m$ with $\mathbf{q}
  \neq \mathbf{q}'$ and let $\epsilon, \epsilon' > 0$ be small.  Then,
  \begin{equation*}
    \mu(B''(\mathbf{q},\epsilon) \cap B''(\mathbf{q}',\epsilon')) \ll 
    \epsilon^n \epsilon'^n.
  \end{equation*}
\end{lem}

\begin{proof}
  First,
  \begin{equation*}
    B''(\mathbf{q},\epsilon) \cap B''(\mathbf{q}',\epsilon') \subseteq
    B(\mathbf{q},\epsilon) \cap B(\mathbf{q}',\epsilon').
  \end{equation*}
  Hence, if $\mathbf{q}$ and $\mathbf{q}'$ are linearly independent, 
  \begin{equation*}
    \mu(B''(\mathbf{q},\epsilon) \cap B''(\mathbf{q}',\epsilon')) \leq 
    \mu(B(\mathbf{q},\epsilon) \cap B(\mathbf{q}',\epsilon')) \asymp
    \epsilon^n {\epsilon'}^n,
  \end{equation*}
  by Lemma \ref{lem:independence} and Lemma \ref{lem:measure}.  Hence,
  we need only show the result when the vectors are linearly
  dependent.  In this case, the result will follow if we can show that
  no overlap between the sets is `too large'.
  
  Let $\hat{\mathbf{q}}$ be a primitive vector in the direction of
  $\mathbf{q}$, \emph{i.e.}, a polynomial vector such that the
  greatest common denominator of the coordinates is of absolute value
  $1$ and such that there are $\lambda, \lambda' \in \mathbb{F}[X]$
  with $\mathbf{q} = \lambda \hat{\mathbf{q}}$ and $\mathbf{q}' =
  \lambda' \hat{\mathbf{q}}$. Suppose without loss of generality that
  $\abs{\lambda} \geq \abs{\lambda'}$, and let $\mathbf{p},
  \mathbf{p}' \in \mathbb{F}[X]^n$ be fixed so that the coordinates of
  $\mathbf{p}$ (respectively those of $\mathbf{p}'$) have greatest
  common denominator of absolute value $1$ with $\lambda$
  (respectively with $\lambda'$). We define affine subspaces $H, H'$
  by
  \begin{alignat*}{2}
    H &= H(\mathbf{q}, \mathbf{p}) = \left\{A \in U : \mathbf{q} A =
      \mathbf{p}\right\}, \\
    H' &= H(\mathbf{q}', \mathbf{p}') = \left\{A' \in
      U : \mathbf{q}' A' = \mathbf{p}'\right\}.
  \end{alignat*}
  It is clear that $B''(\mathbf{q}, \epsilon)$ is the union over
  $\epsilon$-neighbourhoods of all such $H(\mathbf{q}, \mathbf{p})$
  with $\infabs{\mathbf{p}} \leq \infabs{\mathbf{q}}$, where
  $\mathbf{p}$ satisfies the co-primality condition, and similarly for
  $B''(\mathbf{q}', \epsilon')$
  
  To calculate the measure of the intersection, we first note that we
  may disregard contributions from components around such $H$ and $H'$
  if these are separated by a quantity $\geq
  \min\{\epsilon/\infabs{\mathbf{q}},
  \epsilon'/\infabs{\mathbf{q}'}\}$ outside of a set of measure zero.
  We find an arithmetic condition implying this.
  
  Let $A \in H$, $A' \in H'$. On taking the defining equations, we see
  that
  \begin{equation*}
    \lambda \lambda' \hat{\mathbf{q}} (A - A')  = \lambda' \mathbf{p} -
    \lambda \mathbf{p}'. 
  \end{equation*}
  Consider the absolute value of each coordinate of the vector
  $\lambda \lambda' \hat{\mathbf{q}}(A - A')$ in turn . This is some
  expression of the form
  \begin{equation}
    \label{eq:11}
    \abs{\lambda} \abs{\lambda'} \abs{(a_{1j} - a'_{1j}) q_1 + \cdots
      + (a_{mj} - a'_{mj}) q_m}.
  \end{equation}
  If the absolute values of the individual summands have a unique
  maximum, this is the absolute value of the sum by \eqref{eq:27}. We
  show that these absolute values are different outside of a set of
  measure zero.

  Suppose that two summands are of equal absolute value,
  \begin{equation*}
    \abs{(a_{11} - a'_{11}) q_1} = \abs{(a_{12} - a'_{12}) q_2},
  \end{equation*}
  say. Then there is an $\alpha \in \mathbb{F}$ such that 
  \begin{equation*}
    \alpha (a_{11} - a'_{11}) q_1 = (a_{12} - a'_{12}) q_2,
  \end{equation*}
  so that 
  \begin{equation*}
    a_{11} - a_{12}\dfrac{q_2}{\alpha q_1} = a'_{11} - a'_{12}
    \dfrac{q_2}{\alpha q_1}.
  \end{equation*}
  This shows that in order for this to happen, $(a_{11}, a_{12})$ and
  $(a'_{11}, a'_{12})$ must lie on some `line' whose `slope' is a
  rational function. Each such `line' defines an $(nm-1)$-dimensional
  affine space in $\mathbb{F}((X^{-1}))^{mn}$ and so a set of measure
  zero. There are only countably many of these, as there are only
  countably many rational functions over $\mathbb{F}$. Finally, there
  are only finitely many ways in which two summands in \eqref{eq:11}
  can be of equal absolute value. Hence, the exceptional set is of
  measure zero, and we may disregard matrices falling within this set,
  $\mathcal{E}$ say.

  Suppose now that $A, A' \in U \setminus \mathcal{E}$. Then, for all
  $j \in \{1, \dots, n\}$, 
  \begin{equation*}
    \abs{\lambda} \abs{\lambda'} \max_{1 \leq i \leq m}
    \left\{\abs{\hat{\mathbf{q}}_i} \abs{a_{ij} - a_{ij}'}\right\} = 
    \abs{\lambda' p_j - \lambda p_j'}.
  \end{equation*}
  Suppose that for some $j_0 \in \{1, \dots, n\}$, 
  \begin{equation*}
    \abs{\lambda' p_j - \lambda p_j'} \geq \epsilon \abs{\lambda'}.
  \end{equation*}
  It then follows that 
  \begin{multline*}
    \infabs{A - A'} = \max_{\substack{1 \leq i \leq m // 1 \leq j \leq
        n}} \left\{\abs{a_{ij} - a'_{ij}}\right\} \geq \max_{1 \leq i
      \leq m} \left\{\abs{a_{ij_0} - a'_{ij_0}}\right\}\\
    \geq \dfrac{\epsilon}{\abs{\lambda} \infabs{\hat{\mathbf{q}}}}
    \geq \dfrac{\epsilon}{\infabs{\mathbf{q}}} \geq
    \min\left\{\dfrac{\epsilon}{\infabs{\mathbf{q}}},
      \dfrac{\epsilon'}{\infabs{\mathbf{q}'}} \right\},
  \end{multline*}
  which is what was to be shown. By a simple counting argument, it
  then follows that
  \begin{equation}
    \label{eq:10}
    \mu(B''(\mathbf{q},\epsilon) \cap B''(\mathbf{q},\epsilon')) \ll
    \min\left\{\dfrac{\epsilon^n}{\infabs{\mathbf{q}}^n},
      \dfrac{{\epsilon'}^n}{\infabs{\mathbf{q'}}^n} \right\}
    N(\mathbf{q}, \mathbf{q}'), 
  \end{equation}
  where 
  \begin{equation*}
    \begin{split}
      N(\mathbf{q}, \mathbf{q}') = \# \big\{\mathbf{p}, \mathbf{p}'
      \in \mathbb{F}[X]^m :& \infabs{\mathbf{p}} \leq
      \infabs{\mathbf{q}}, \infabs{\mathbf{p}} \leq
      \infabs{\mathbf{q}},\\ 
      & (p_j, \lambda) =
      (p'_j, \lambda') = 1 \\
      &\text{and } \abs{\lambda' p_j - \lambda p_j'} < \epsilon
      \abs{\lambda'} \text{ for all } j = 1, \dots, n \big\}
    \end{split}
  \end{equation*}
  Arguing exactly as in \cite[pp.~158--159]{MR2019008} for each
  coordinate, we find that
  \begin{equation}
    \label{eq:29}
    N(\mathbf{q}, \mathbf{q}') \leq \left(\infabs{\mathbf{q}'}
    \epsilon + \infabs{\mathbf{q}} \epsilon'\right)^n. 
  \end{equation}
  Inserting \eqref{eq:29} into \eqref{eq:10}, we obtain the desired
  inequality. 
\end{proof}

Finally, we will construct a subset $\tilde{S} \subseteq S$ which
contains a lot of elements. First, define for $N \in \mathbb{N}$,
\begin{equation*}
  S_N = \left\{\mathbf{q} \in S : k^N \leq \infabs{\mathbf{q}} <
    k^{N+1}\right\}.
\end{equation*}
We will construct a subset $\tilde{S} \subseteq S$ and will denote by
$\tilde{S}_N$ the $N$-th $k$-adic block as above.

\begin{lem}
  \label{lem:large_blocks}
  Let $\delta > 0$. For any $N_0 \in \mathbb{N}$, there is an $N \geq
  N_0$ such that 
  \begin{equation*}
    \# S_N \geq k^{N(v(S)-\delta)}.
  \end{equation*}
\end{lem}

\begin{proof}
  Suppose to the contrary that for for some $N_0 \in \mathbb{N}$, 
  \begin{equation*}
    \# S_N < k^{N(v(S)-\delta)},
  \end{equation*}
  whenever $N \geq N_0$. Then, for any fixed $\epsilon < \delta$,
  \begin{multline*}
    \sum_{\substack{\mathbf{q} \in S \\ \infabs{\mathbf{q}} \geq
        k^{N_0}}} \infabs{\mathbf{q}}^{-v(S)+\epsilon} =
    \sum_{N=N_0}^\infty \sum_{\mathbf{q} \in S_N}
    \infabs{\mathbf{q}}^{-v(S)+\epsilon} \ll \sum_{N= N_0}^\infty
    k^{-N(v(S) - \epsilon)} \# S_N \\
    < \sum_{N= N_0}^\infty k^{-N(v(S) - \epsilon) + N(v(S)-\delta)} =
    \sum_{N= N_0}^\infty \left(k^{\epsilon-\delta}\right)^N < \infty.  
  \end{multline*}
  On the other hand, by definition, 
  \begin{equation*}
    \sum_{\substack{\mathbf{q} \in S \\ \infabs{\mathbf{q}} \geq
        k^{N_0}}} \infabs{\mathbf{q}}^{-v(S)+\epsilon} = \infty,
  \end{equation*}
  for any $\epsilon> 0$.
\end{proof}

\subsection{Probabilistic lemmas}
\label{sec:probabilistic-lemmas}

We will define random variables which will be used to construct a
ubiquitous system. First, suppose that $m = 1$, let $\delta > 0$ be
arbitrary and let $\{N_t\}_{t=1}^\infty$ be an increasing sequence of
integers such that Lemma \ref{lem:large_blocks} holds for each $N_t$.
Finally, define the function
\begin{equation}
  \label{eq:12}
  \rho(k^N) = k^{-N(v(S) - \delta)/n}\log N.    
\end{equation}
For each $t \in \mathbb{N}$, we define random variables on $U$,
\begin{equation*}
  \nu_t(A) = \sum_{q \in S_{N_r}} \chi_{B'(q, \rho(k^{N_t}))}(A),
\end{equation*}
where $\chi_E$ denotes the characteristic function of the set $E$.
Clearly, these variables count the number of sets $B'(q,
\rho(k^{N_t}))$ with $q \in S_{N_t}$ which contain a given element $A
\in U$. Consequently,
\begin{equation*}
  \nu_t^{-1}(0) = U \setminus \bigcup_{q \in S_{N_t}} B'(q,
  \rho(k^{N_t})). 
\end{equation*}
We calculate the first and second moment of these random variables. 

First, let $r \in \mathbb{Z}$ be such that $k^r \leq \rho(k^{N_t}) <
k^{r+1}$. This implies that
\begin{equation*}
  B'(q,\rho(k^{N_t})) = B'(q, k^r).
\end{equation*}
We easily estimate the mean value by integrating,
\begin{equation}
  \label{eq:5}
  \begin{split}
    \mathbb{E}(\nu_t) &= \int_U \sum_{\mathbf{q} \in S_{N_t}}
    \chi_{B'(q,\rho(k^N))}(X) d\mu(X) \\
    &= \sum_{q \in S_{N_t}} \int_U \chi_{B'(q,k^r)}(X) d\mu(X) 
    \asymp \sum_{q \in S_{N_t}} k^{rn} \asymp \rho(k^{N_t})^n \#
    S_{N_t}, 
  \end{split}
\end{equation}
by Lemma \ref{lem:1D-measure} and Lemma \ref{lem:phi-function}.

We now estimate the second moment using the pairwise
quasi-in\-de\-pen\-dence property from Lemma
\ref{lem:1D-intersection}, the measure estimate from Lemma
\ref{lem:1D-measure} and Lemma~\ref{lem:phi-function}.
\begin{equation}
  \label{eq:7}
  \begin{split}
    \mathbb{E}\left(\nu_t^2\right) &= \int_U \sum_{q \in
      S_{N_t}} \chi_{B'(q,\rho(k^{N_t}))}(X)^2 d\mu(X)\\
    &\quad\quad + \int_U \sum_{\substack{q, q' \in S_{N_t} \\ q \neq q'}}
    \chi_{B'(q,\rho(k^{N_t}))}(X) \chi_{B'(q',\rho(k^{N_t}))}(X) d\mu(X) 
    \\
    &\ll \mathbb{E}(\nu_t) + \sum_{\substack{q, q' \in S_{N_t} \\
        q \neq q'}} \rho(k^{N_t})^n \rho(k^{N_t})^n \leq
    \mathbb{E}(\nu_t) + \mathbb{E}(\nu_t)^2.
  \end{split}
\end{equation}
By \eqref{eq:5} and \eqref{eq:7}, the variance $\sigma_t^2$ of $\nu_t$
satisfies
\begin{equation}
  \label{eq:8}
  \sigma_t^2 = \mathbb{E}\left(\nu_t^2\right) -
  \mathbb{E}(\nu_t)^2 \leq \mathbb{E}(\nu_t).
\end{equation}
This has the following consequence.

\begin{lem}
  \label{lem:pre-ubiquity}
  With $\nu_t$ and $\rho(N)$ as above,
  \begin{equation*}
    \mu(\nu_t^{-1}(0)) \leq \dfrac{1}{\mathbb{E}(\nu_t)} \rightarrow 0,
  \end{equation*}
  as $t$ tends to infinity.
\end{lem}

\begin{proof}
  The proof is easy. We use an alternative characterisation of the
  variance: 
  \begin{multline*}
    \sigma^2_{t} = \int_{U} (v_{t}(A) - \mathbb{E}(\nu_t))^2 d\mu(A)
    \\ 
    \geq \int_{\nu_{t}^{-1}(0)} (v_{t}(A) - \mathbb{E}(\nu_t))^2
    d\mu(A) = \mathbb{E}(\nu_t)^2 \mu\left(\nu_{t}^{-1}(0)\right).
  \end{multline*}
  Hence, 
  \begin{equation*}
    \mu\left(\nu_{t}^{-1}(0)\right) \leq
    \dfrac{\sigma^2_{t}}{\mathbb{E}(\nu_t)^2} \leq
    \dfrac{1}{\mathbb{E}(\nu_t)},
  \end{equation*}
  by \eqref{eq:8}. Finally, note that by \eqref{eq:5} and Lemma
  \ref{lem:large_blocks}, 
  \begin{equation*}
    \dfrac{1}{\mathbb{E}(\nu_t)} \asymp \dfrac{1}{\rho(k^{N_t})^n \#
      S_{N_t}} \ll \dfrac{1}{(\log N_t)^n}.
  \end{equation*}
  As $N_t$ is increasing, this completes the proof.
\end{proof}

We now proceed to discuss the case when $m \geq 2$.  As before, we let
$\delta > 0$ be arbitrary and let $\{N_t\}_{t=1}^\infty$ be an
increasing sequence of integers such that Lemma \ref{lem:large_blocks}
holds for each $N_t$. For each $t \in \mathbb{N}$, we define random
variables on $U$,
\begin{equation*}
  \tilde{\nu}_t(A) = \sum_{q \in \tilde{S}_{N_r}} \chi_{B''(q,
    \rho(k^{N_t}))}(A). 
\end{equation*}
These may be interpreted as $\nu_t$ above.

Calculating the first and second moments may be done exactly as in
\eqref{eq:5} and \eqref{eq:7} by using the measure estimates of Lemma
\ref{lem:big_m_measure} and Lemma \ref{lem:big_m_intersection}, so
that
\begin{equation*}
  \mathbb{E}(\tilde{\nu}_t) \asymp \rho(k^{N_t})^n \# \tilde{S}_{N_t}, 
\end{equation*}
and
\begin{equation*}
  \mathbb{E}\left(\tilde{\nu}_t^2\right) \leq
  \mathbb{E}(\tilde{\nu}_t) + \mathbb{E}(\tilde{\nu}_t)^2. 
\end{equation*}
We recognise this as the main ingredients in the proof of Lemma
\ref{lem:pre-ubiquity}, so that we immediately obtain the analogous
version for $m \geq 2$.
\begin{lem}
  \label{lem:big_m_pre-ubiquity}
  With $\tilde{\nu}_t$ and $\rho(N)$ as above,
  \begin{equation*}
    \mu(\tilde{\nu}_t^{-1}(0)) \leq
    \dfrac{1}{\mathbb{E}(\tilde{\nu}_t)} \rightarrow 0,
  \end{equation*}
  as $t$ tends to infinity.
\end{lem}

\subsection{Completing the proof}
\label{sec:completing-proof}

We first show \eqref{item:1}, so suppose that $nv < v(S)$. The following
result is classical: If $\{A_i\}$ is an infinite sequence of events in
some probability space with probability measure $\mathbb{P}$, such
that $\sum \mathbb{P}(A_i) = \infty$ and $A = \bigcap_{N=1}^\infty
\bigcup_{i=N}^\infty A_i$, then
\begin{equation}
  \label{eq:13}
  \mathbb{P}(A) \geq \limsup_{N \rightarrow \infty}
  \dfrac{\left(\sum_{i=1}^N \mathbb{P}(A_i)\right)^2}{\sum_{i=1}^N
    \sum_{j=1}^N \mathbb{P}(A_i \cap A_j)}.
\end{equation}
We will apply this to the probability space $U$ equipped with the Haar
measure and with the events $B'(q, \abs{q}^{-v})$ (resp.
$B''(\mathbf{q}, \infabs{\mathbf{q}}^{-v})$).

First, let $m=1$ and fix $\delta < v(S) - nv$. By Lemma
\ref{lem:1D-measure} and Lemma \ref{lem:large_blocks},
\begin{equation*}
  \sum_{q \in \mathbb{F}[X]} \mu(B'(q, \abs{q}^{-v})) \gg
  \sum_{t=1}^\infty \sum_{q \in S_{N_t}} \abs{q}^{-nv} \asymp
  \sum_{t=1}^\infty k^{N_t (v(S) - \delta -nv)} = \infty.
\end{equation*}
Hence, \eqref{eq:13} applies. Using Lemma \ref{lem:1D-intersection},
we find that
\begin{equation*}
  \mu(\mathcal{W}^*_S(1,n;\abs{ \cdot }^{-v})) \geq c > 0
\end{equation*}
where $c$ comes from the right hand side of \eqref{eq:13}. Making the
obvious modifications when $m \geq 2$, we find for all $m,n$, 
\begin{equation}
  \label{eq:19}
  \mu(\mathcal{W}^*_S(m,n;\infabs{ \cdot }^{-v})) \geq c > 0
\end{equation}

We now apply an inflation argument due to Cassels \cite{MR0087708} to
show that the measure must be full. This is in complete analogy with
the proof of \cite[Theorem 3]{kristensen03}. In the above argument,
we could just as easily have shown that
$\mu(\mathcal{W}^*_S(m,n;\eta(\cdot) \infabs{ \cdot }^{-v})) \geq c >
0$, where $\eta: \mathbb{F}[X]^m \rightarrow (0,1]$ is a function
depending only on $\infabs{\mathbf{q}}$ which decreases to zero as
$\infabs{\mathbf{q}}$ increases, such that $\lambda(\eta(\cdot)
\infabs{ \cdot }^{-v}) = v$.  Hence, can find a point of metric
density $A_0 \in \mathcal{W}^*_S(m,n;\eta(\cdot) \infabs{ \cdot
}^{-v})$. Let $\epsilon > 0$.  We may find an $r_0 \in \natnum$ with
\begin{equation*}
  \mu \left(\mathcal{W}^*_S(m,n;\eta(\cdot) \infabs{ \cdot }^{-v}) \cap
    B(A_0, k^{-r_0})\right) \geq k^{-mn r_0}-\dfrac{\epsilon}{k^{mn
      r_0}}.
\end{equation*}
We scale the set by $X^{r_0}$ to obtain
\begin{equation*}
  \mu\left( X^{r_0} \left(\mathcal{W}^*_S(m,n;\eta(\cdot) \infabs{ \cdot 
      }^{-v}) \cap B(A_0, k^{-r_0})\right)\right) \geq 1- \epsilon, 
\end{equation*}
since $\mu(X^{r_0} B(c,r)) = k^{mn r_0} \mu(B(c,r))$ for any ball
$B(c,r)$.  

It follows that for any $\epsilon > 0$, we may find a set $S \subseteq
U$ of measure $\mu(S) \geq 1-\epsilon$, such that any $A \in S$ may be
written on the form
\begin{equation}
  \label{eq:24}
  A = X^{r_0} A' + P, \quad A' \in \mathcal{W}^*_S(m,n;\eta(\cdot)
  \infabs{ \cdot }^{-v}), P \in \mathbb{F}[X]^{mn}.
\end{equation}
On considering the distance to nearest polynomial vector for such
elements, we find that $S \subseteq \mathcal{W}^*_S(m,n;k^{r_0}
\eta(\cdot) \infabs{ \cdot }^{-v})$. Considering for fixed $A \in
\mathcal{W}^*_S(m,n;k^{r_0} \eta(\cdot) \infabs{ \cdot }^{-v})$, the
(infinitely many) $\mathbf{q}$ for which $\eta(\mathbf{q}) < k^{-r_0}$
and for which $A<k^{r_0} \eta(\mathbf{q}) \infabs{\mathbf{q}}^{-v}$
shows that 
\begin{equation*}
  S \subseteq \mathcal{W}^*_S(m,n;k^{r_0} \eta(\cdot) \infabs{ \cdot
  }^{-v}) \subseteq \mathcal{W}^*_S(m,n;\infabs{ \cdot }^{-v}).
\end{equation*}
This finishes the proof.

We now prove the lower bound on the dimension in \eqref{item:2}, which
together with the upper bound found in \S\ref{sec:an-upper-bound} will
complete the proof of the theorem. We will use the above estimates to
construct a ubiquitous system of sets. 

We define a function $\tilde{\rho}(k^N) = \rho(k^N) k^{-N+1}$, where
$\rho$ is the function defined in \eqref{eq:12}, and sets
\begin{equation*}
  \tilde{B}(\mathbf{q};\epsilon) = \left\{A \in U : \dist\left(A,
      \bigcup_{\mathbf{p} \in \mathbb{F}[X]^n} H(\mathbf{q},
      \mathbf{p})\right) < \epsilon  \right\},
\end{equation*}
where $H(\mathbf{q}, \mathbf{p}) = \{ A \in U :
\mathbf{q}A=\mathbf{p}\}$.  We shall prove that when $\mathbf{q} \in
S_{N_t}$, then
\begin{equation}
  \label{eq:20}
  B(\mathbf{q};\rho(k^{N_t})) \subseteq
  \tilde{B}(\mathbf{q};\tilde{\rho}(k^{N_t})). 
\end{equation}
But this easily follows as clearly, 
\begin{equation*}
  \infabs{\mathbf{q}} \dist(A, H(\mathbf{q},\mathbf{p})) \leq
  \infabs{\mathbf{q}A - \mathbf{p}}.
\end{equation*}
Choosing $\mathbf{p}$ so that $\infabs{\mathbf{q} A-\mathbf{p}} <
\rho(k^{N_t})$, and noting that $k^{N_t-1} \leq \infabs{\mathbf{q}}$
as $\mathbf{q} \in S_{N_t}$, we have shown \eqref{eq:20}.

We now claim that the system $(\bigcup_{p \in \mathbb{F}[X]^n}
H(\mathbf{q}, \mathbf{p}), \infabs{\mathbf{q}})$, where $\mathbf{q}$
runs over $\bigcup_{t=1}^\infty S_{N_t}$ is ubiquitous with respect to
$\tilde{\rho}(N)$, \emph{i.e.},
\begin{equation*}
  \lim_{t \rightarrow \infty} \mu\left(U \setminus \bigcup_{1 \leq
      \infabs{\mathbf{q}} \leq k^{N_t}}
    \tilde{B}(\mathbf{q};\tilde{\rho}(k^{N_t}))\right) = 0.
\end{equation*}
But this follows from \eqref{eq:20}, as
\begin{alignat*}{2}
  U \setminus \bigcup_{1 \leq\infabs{\mathbf{q}} \leq k^{N_t}}
  \tilde{B}(\mathbf{q};\tilde{\rho}(k^{N_t})) &\subseteq U \setminus  
  \bigcup_{\mathbf{q} \in S_{N_t}} B\left(\mathbf{q},
    \rho(k^{N_t})\right)\\
  &\subseteq 
  \begin{cases}
    U \setminus \bigcup_{\mathbf{q} \in S_{N_t}} B'\left(\mathbf{q}, 
      \rho(k^{N_t})\right) & \text{for $m = 1$},\\
    U \setminus \bigcup_{\mathbf{q} \in S_{N_t}} B''\left(\mathbf{q}, 
      \rho(k^{N_t})\right) & \text{for $m \geq 2$}.\\
  \end{cases}
\end{alignat*}
By Lemma \ref{lem:pre-ubiquity} and Lemma
\ref{lem:big_m_pre-ubiquity}, the measure of the right hand side tends
to zero as $t$ tends to infinity.

Now using the above, \cite[Lemma 6]{kristensen03} implies that 
\begin{alignat*}{2}
  \hdim(\mathcal{W}^*_S&(m,n;\psi))  \geq n(m-1) + n
  \limsup_{t \rightarrow \infty} \dfrac{\log_k
    \tilde{\rho} (k^{N_t})}{\log_k(k^{N_t(-v-1)})} \\
  & = n(m-1) +\limsup_{t \rightarrow \infty} \dfrac{-N_t (v(S) + n
    - \delta) -1 +\log_k \log N_t}{N_t(-v-1)} \\
  &= n(m-1) + \dfrac{v(S) + n - \delta}{v+1}.
\end{alignat*}
As $\delta$ was arbitrary, this completes the proof of Theorem
\ref{thm:Main}. \qed

\section{Proof of Theorem \ref{thm:Main2}}
\label{sec:proof-theorem-2}

We now prove Theorem \ref{thm:Main2}. In order to accomplish this, we
essentially use Rynne's method \cite{MR1633797}. First, note that if
$\psi(\mathbf{q}) > 1$ for infinitely many $\mathbf{q} \in
\mathbb{F}[X]^m$, the theorem is trivially true. Hence, there is no
loss of generality in supposing that $\psi(\mathbf{q}) \leq 1$ for all
$\mathbf{q} \in \mathbb{F}[X]^m$, as we may ignore the finitely many
cases for which this happens. We will prove four lemmas, which will
imply the theorem. First, we show that the upper bound on the
dimension is the right one.

\begin{lem}
  \label{lem:upper_bound_DS}
  With $n,m,\psi$ and $\eta(\psi)$ as above, 
  \begin{equation*}
    \hdim(\mathcal{W}(m,n;\psi)) \leq n(m-1) + \min\{\eta(\psi),
    n\}.
  \end{equation*}
\end{lem}

\begin{proof}
  This is identical to the covering argument of
  \S\ref{sec:an-upper-bound}, where we replace the balls of the cover
  by balls of radius $\psi(\mathbf{q})/\infabs{\mathbf{q}}$. This does
  not affect the number of balls needed, and we do not repeat the
  argument here.
\end{proof}

In order to show that the lower bound hold, we define a number of
quantities to be used in the proof. Let $\psi$ be a function as in the
statement of Theorem \ref{thm:Main2}, let $N \in \mathbb{N}$ and let
$v > 0$. We define
\begin{alignat*}{2}
  C(N, v; \psi) &= \#\{\mathbf{q} \in \mathbb{F}[X]^m :
  \infabs{\mathbf{q}} \leq N, \psi(\mathbf{q}) \geq
  \infabs{\mathbf{q}}^{-v}\}, \\
  \gamma(v; \psi)&= \sup\{\gamma \in \mathbb{R} : \limsup_{N
    \rightarrow \infty} C(N, v;\psi)N^{-\gamma} > 0 \}, \\
  & \quad\quad \text{ defined whenever } \lim_{N
    \rightarrow \infty} C(N, v; \psi) = \infty, \\
  \delta(v;\psi) &=
  \begin{cases}
    \dfrac{n + \gamma(v; \psi)}{v+1} & \text{if } \lim_{N
      \rightarrow \infty} C(N, v; \psi) = \infty, \\
    0 & \text{otherwise,}
  \end{cases}\\
  \delta(\psi) &= \sup_{v \geq 0} \delta(v;\psi).
\end{alignat*}
For arbitrary infinite subsets $S \subseteq \mathbb{F}[X]^m$, we will
need an additional two definitions,
\begin{alignat*}{2}
  C(N; S) &= \# \{\mathbf{q} \in S : \infabs{\mathbf{q}} \leq N\}, \\
  \gamma(S) &= \sup\{\gamma \in \mathbb{R} : \limsup_{N \rightarrow
    \infty} C(N;S)N^{-\gamma} > 0 \}.
\end{alignat*}

We first show that the exponent $\gamma(S)$ coinsides with the
exponent $v(S)$ of Theorem \ref{thm:Main}.

\begin{lem}
  \label{lem:exponents}
  Let $S \subseteq \mathbb{F}[X]^m$ with $\#S = \infty$. Then, $v(S) =
  \gamma(S)$.
\end{lem}

\begin{proof}
  Let $\epsilon > 0$ be fixed. Choose $N_0$ so that $C(N;S) \leq
  N^{\gamma(S)+\epsilon}$ for $N \geq N_0$ and let $v > \gamma(S) +
  \epsilon$. Then, 
  \begin{multline*}
    \sum_{\mathbf{q} \in Q} \infabs{\mathbf{q}}^{-v} =
    \sum_{r=1}^\infty \sum_{k^{r-1} \leq \infabs{\mathbf{q}} < k^r} 
    \infabs{\mathbf{q}}^{-v} \ll \sum_{r=1}^\infty C(k^r;S) k^{-rv} \\ 
    \ll \sum_{r=1}^\infty k^{r(\gamma(S) + \epsilon - v)} < \infty.
  \end{multline*}
  Hence, $v(S) \leq \gamma(S)$.

  For the reverse inequality, let $\gamma < \gamma(S)$ and choose a
  sequence $N_r$ with 
  \begin{equation*}
    C(N_r;S) N_r^{-\gamma} > \epsilon \quad \text{and} \quad
    C(N_{r+1};S) > 2 C(N_r;S)
  \end{equation*}
  for all $r \in \mathbb{N}$. Clearly this is possible, as $\# S =
  \infty$. Now, 
  \begin{multline*}
    \sum_{\mathbf{q} \in Q} \infabs{\mathbf{q}}^{-\gamma} \geq
    \sum_{r=1}^\infty \sum_{N_r \leq \infabs{\mathbf{q}} < N_{r+1}}
    \infabs{\mathbf{q}}^{-\gamma} \geq \sum_{r=1}^\infty (C(N_{r+1};S)
    - C(N_r;S)) N_{r+1}^{-\gamma} \\
    > \sum_{r=1}^\infty \tfrac{1}{2} C(N_{r+1};S) N_{r+1}^{-\gamma} =
    \tfrac{1}{2} \sum_{r=1}^\infty \epsilon = \infty.
  \end{multline*}
  Hence, $\gamma < v(S)$, so that $\gamma(S) \leq v(S)$.
\end{proof}

\begin{lem}
  \label{lem:lower_bound_DS}
  With $n,m,\psi$ and $\delta(\psi)$ as above, 
  \begin{equation*}
    \hdim(\mathcal{W}(m,n;\psi)) \geq n(m-1) + \min\{\delta(\psi),
    n\}.
  \end{equation*}
\end{lem}

\begin{proof}
  We will deduce this lemma from Theorem \ref{thm:Main}. The set
  $\mathcal{W}(m,n;\psi)$ contains $n(m-1)$-dimensional affine spaces,
  and so the dimension is at least $n(m-1)$. Furthermore, it is at
  most $mn$, so we may suppose that $0 < \delta(\psi) \leq n$.

  Let $0 < \epsilon < \delta(\psi)$ and fix $v_0 \geq 0$ so that $\delta(v_0;\psi) >
  \delta(\psi) -\epsilon > 0$. Define a set 
  \begin{equation*}
    S = \{\mathbf{q} \in \mathbb{F}[X]^m : \psi(\mathbf{q}) \geq
    \infabs{\mathbf{q}}^{-v_0}\}.
  \end{equation*}
  This set is infinite by choice of $v_0$. Also, 
  \begin{equation}
    \label{eq:21}
    \mathcal{W}_S(m,n; \infabs{\cdot}^{-v_0}) \subseteq
    \mathcal{W}(m,n;\psi).
  \end{equation}
  Finally, $C(N, v_0; \psi) = C(N;S) \rightarrow \infty$ as $N$ tends
  to infinity, so that 
  \begin{equation}
    \label{eq:22}
    \gamma(v_0; \psi) = \gamma(S) = v(S),
  \end{equation}
  by Lemma \ref{lem:exponents}. 

  We apply Theorem \ref{thm:Main} to $\mathcal{W}_S(m,n;
  \infabs{\cdot}^{-v_0})$. In view of \eqref{eq:21} and \eqref{eq:22},  
  \begin{multline*}
    \hdim \left(\mathcal{W}(m,n;\psi)\right) \geq n(m-1) + \dfrac{n +
      v(S)}{1+v_0} \\
    = n(m-1) + \delta(v_0; \psi) > n(m-1)+\delta(\psi)
    -\epsilon. 
  \end{multline*}
  As $\epsilon$ was arbitrary, this completes the proof.
\end{proof}

\begin{lem}
  \label{lem:compare_exponents}
  With $n,m,\psi,\delta(\psi)$  and $\eta(\psi)$ as above, 
  \begin{equation*}
    \min\{\delta(\psi), n\} \geq \min\{\eta(\psi), n\}.
  \end{equation*}
\end{lem}

\begin{proof}
  Let $\eta > \delta(\psi)$. We will show that $\eta(\psi) \leq \eta$,
  so that we must have $\eta(\psi) \leq \delta(\psi)$. In order to
  accomplish this, we will split the series
  \begin{equation*}
    \sum_{\mathbf{q} \in \mathbb{F}[X]^m} \infabs{\mathbf{q}}^n \left( 
      \dfrac{\psi(\mathbf{q})}{\infabs{\mathbf{q}}} \right)^\eta
  \end{equation*}
  up into finitely many components, each of which converge. This will
  imply the result.

  First, let $\mu = (m+n)/\eta$. Consider the set
  \begin{equation*}
    S' = \left\{\mathbf{q} \in \mathbb{F}[X]^m : \psi(\mathbf{q}) <
      \infabs{\mathbf{q}}^{-\mu}\right\}.
  \end{equation*}
  Then, by \cite[equation (1.4)]{kristensen03},
  \begin{equation}
    \label{eq:23}
    \sum_{\mathbf{q} \in S'} \infabs{\mathbf{q}}^n \left( 
      \dfrac{\psi(\mathbf{q})}{\infabs{\mathbf{q}}} \right)^\eta \leq
    \sum_{\mathbf{q} \in \mathbb{F}[X]^m} \infabs{\mathbf{q}}^{n-(\mu
      +1)\eta} \ll \sum_{r=1}^\infty k^{r(m+n-(\mu+1)\eta)} < \infty. 
  \end{equation}
  
  Now, fix a $\theta \in (0, 1-\delta(\psi)/\eta)$. Let $N_0$ be such
  that $C(N;S) \leq N^{\gamma(v;\psi) + \epsilon}$ for $N \geq N_0$.
  For any $v > 0$, we define sets
  \begin{equation*}
    S(v,\theta) = \left\{\mathbf{q} \in \mathbb{F}[X]^m :
      \infabs{\mathbf{q}}^{-v} \leq \psi(\mathbf{q}) \leq
      \infabs{\mathbf{q}}^{-v+\theta}\right\}. 
  \end{equation*}
  Clearly, since we can cover the set $[0, \mu]$ by intervals of the
  form $[v-\theta, v]$, there is a finite set $V$, such that
  \begin{equation}
    \label{eq:25}
    \mathbb{F}[X]^m = S' \cup \bigcup_{v \in V} S(v,\theta),
  \end{equation}
  since $\psi$ is assumed to be bounded above by $1$.  Choose such a
  finite set $V$ and let $0 < \epsilon < (\eta-\delta(\psi)) \min_{v
    \in V} \{v\}$. For any $v \in V$,
  \begin{alignat*}{2}
    \sum_{\mathbf{q} \in S(v, \theta)} \infabs{\mathbf{q}}^n \left( 
      \dfrac{\psi(\mathbf{q})}{\infabs{\mathbf{q}}} \right)^\eta &\leq  
    \sum_{\mathbf{q} \in S(v,\theta)} \infabs{\mathbf{q}}^{n-(v
      +1 - \theta)\eta} \\
    &\ll \sum_{r=1}^\infty \sum_{\substack{k^{r-1}
        \leq \infabs{q} \leq k^r \\ \mathbf{q} \in S(v,\theta)}} k^{r(n-(v
      +1 - \theta)\eta)} \\
    & \ll \sum_{r=1}^\infty k^{r(n-(v +1 - \theta)\eta +
      \gamma(v,\psi) + \epsilon)}.
  \end{alignat*}
  Hence, to show that this series converges, it suffices to show that
  the exponent is negative. But this follows as
  \begin{alignat*}{2}
    n-(v +1 - \theta)\eta + \gamma(v,\psi) + \epsilon 
    &< n - v \eta - \delta(\psi) + \left((v+1)\delta(v,\psi) - n\right)
    + \epsilon\\
    &< (\delta(\psi) - \eta)v + \epsilon
  \end{alignat*}
  which is less than zero by choice of $\epsilon$. Hence,
  \begin{equation}
    \label{eq:26}
    \sum_{\mathbf{q} \in S(v, \theta)} \infabs{\mathbf{q}}^n \left( 
      \dfrac{\psi(\mathbf{q})}{\infabs{\mathbf{q}}} \right)^\eta <
    \infty, 
  \end{equation}
  and the lemma follows by \eqref{eq:23}, \eqref{eq:25} and
  \eqref{eq:26}. 
\end{proof}

We now complete the proof of Theorem \ref{thm:Main2} by noting that
the upper bound follows immediately from Lemma
\ref{lem:upper_bound_DS}. The lower bound follows on inserting the
result of Lemma \ref{lem:compare_exponents} into Lemma
\ref{lem:lower_bound_DS}. \qed

\section{Acknowledgements}
\label{sec:acknowledgements}

I thank Yann Bugeaud for pointing out the paper by Inoue and Nakada
\cite{inoue03:_dioph}. I also thank an anonymous referee for
suggesting substantial improvements on the original
manuscript. Finally, I thank Department of Mathematical Sciences at
the University of Aarhus for their kind hospitality.


\begin{thebibliography}{10}

\bibitem{MR2001h:11091}
V.~I. Bernik and M.~M. Dodson.
\newblock {\em Metric {D}iophantine approximation on manifolds}, volume 137 of
  {\em Cambridge Tracts in Mathematics}.
\newblock Cambridge University Press, Cambridge, 1999.

\bibitem{MR46:7175}
I.~Borosh and A.~S. Fraenkel.
\newblock A generalization of {J}arn\'\i k's theorem on {D}iophantine
  approximations.
\newblock {\em Nederl. Akad. Wetensch. Proc. Ser. A {\bf 75}=Indag. Math.},
  34:193--201, 1972.

\bibitem{MR0087708}
J.~W.~S. Cassels.
\newblock {\em An introduction to {D}iophantine approximation}.
\newblock Cambridge Tracts in Mathematics and Mathematical Physics, No. 45.
  Cambridge University Press, New York, 1957.

\bibitem{MR43:161}
B.~de~Mathan.
\newblock Approximations diophantiennes dans un corps local.
\newblock {\em Bull. Soc. Math. France Suppl. M\'em.}, 21:93, 1970.

\bibitem{MR2001g:11129}
H.~Dickinson and B.~P. Rynne.
\newblock Hausdorff dimension and a generalized form of simultaneous
  {D}iophantine approximation.
\newblock {\em Acta Arith.}, 93(1):21--36, 2000.

\bibitem{MR95d:11092}
M.~M. Dodson.
\newblock Geometric and probabilistic ideas in the metric theory of
  {D}iophantine approximations.
\newblock {\em Uspekhi Mat. Nauk}, 48(5(293)):77--106, 1993.

\bibitem{MR2019008}
K.~Inoue.
\newblock The metric simultaneous {D}iophantine approximations over formal
  power series.
\newblock {\em J. Th\'eor. Nombres Bordeaux}, 15(1):151--161, 2003.
\newblock Les XXII\`emes Journ\'ees Arithmetiques (Lille, 2001).

\bibitem{inoue03:_dioph}
K.~Inoue and H.~Nakada.
\newblock On metric {D}iophantine approximation in positive characteristic.
\newblock {\em Acta Arith.}, 110(3):205--218, 2003.

\bibitem{kristensen03}
S.~Kristensen.
\newblock On well-approximable matrices in a field of formal series.
\newblock {\em Math. Proc. Cambridge Philos. Soc.}, 135(2):255--268, 2003.

\bibitem{MR2001k:11135}
A.~Lasjaunias.
\newblock A survey of {D}iophantine approximation in fields of power series.
\newblock {\em Monatsh. Math.}, 130(3):211--229, 2000.

\bibitem{MR1099767}
A.~D. Pollington and R.~C. Vaughan.
\newblock The {$k$}-dimensional {D}uffin and {S}chaeffer conjecture.
\newblock {\em Mathematika}, 37(2):190--200, 1990.

\bibitem{MR93a:11066}
B.~P. Rynne.
\newblock The {H}ausdorff dimension of certain sets arising from {D}iophantine
  approximation by restricted sequences of integer vectors.
\newblock {\em Acta Arith.}, 61(1):69--81, 1992.

\bibitem{MR1633797}
B.~P. Rynne.
\newblock The {H}ausdorff dimension of sets arising from {D}iophantine
  approximation with a general error function.
\newblock {\em J. Number Theory}, 71(2):166--171, 1998.

\bibitem{MR2001j:11063}
W.~M. Schmidt.
\newblock On continued fractions and {D}iophantine approximation in power
  series fields.
\newblock {\em Acta Arith.}, 95(2):139--166, 2000.

\bibitem{MR548467}
V.~G. Sprind{\v{z}}uk.
\newblock {\em Metric theory of {D}iophantine approximations}.
\newblock V. H. Winston \& Sons, Washington, D.C., 1979.

\end{thebibliography}
\end{document}